\documentclass[a4paper,11pt]{article}
\title{A Unique Perfect Power Decagonal Number}
\author{Philippe Michaud-Rodgers}

\makeatletter
\newcommand\notsotiny{\@setfontsize\notsotiny\@vipt\@viipt}
\makeatother
\usepackage{amsmath} 
\usepackage{amssymb} 
\usepackage{mathtools}
\usepackage{framed}
\usepackage{pifont} 
\usepackage{enumerate} 
\usepackage{wasysym} 
\usepackage{commath}
\usepackage{graphicx}
\usepackage{framed}
\usepackage{amsmath}
\usepackage{amssymb}
\usepackage{amsfonts}
\usepackage{mathrsfs}
\usepackage{mathtools}
\usepackage{url}
\usepackage{array}
\usepackage{amsthm}
\usepackage[english]{babel}
\usepackage[utf8]{inputenc}
\newtheorem{theorem}{Theorem}[section]
\newtheorem*{theorem*}{Theorem}

\theoremstyle{definition}

\theoremstyle{remark}
\newtheorem{remark}[theorem]{Remark}
\newtheorem*{remark*}{Remark}
\usepackage{yfonts}
\usepackage{amsmath,amssymb,tikz-cd}
\usepackage{pdflscape}

\usepackage{etoolbox}
\apptocmd{\sloppy}{\hbadness 10000\relax}{}{}

\usepackage[numbers]{natbib}
\usepackage{cite}
\usepackage{mathtools, nccmath}

\providecommand{\Q}{\mathbb{Q}}
\providecommand{\Z}{\mathbb{Z}}

\usepackage{longtable}
\usepackage{hyperref}

\newcommand{\Addresses}{{
  \bigskip
  \footnotesize

 \textsc{Mathematics Institute, University of Warwick, CV4 7AL, United Kingdom}\par\nopagebreak
  \textit{E-mail address}: \texttt{p.rodgers@warwick.ac.uk}
}}

\begin{document}

\maketitle
 
\begin{abstract}
Let  $\mathcal{P}_s(n)$ denote the $n$th $s$-gonal number. We consider the equation \[\mathcal{P}_s(n) = y^m, \] for integers $n,s,y,$ and $m$. All solutions to this equation are known for $m>2$ and $s \in \{3,5,6,8,20 \}$. We consider the case $s=10$, that of decagonal numbers. Using a descent argument and the modular method, we prove that the only decagonal number $>1$ expressible as a perfect $m$th power with $m>1$ is $\mathcal{P}_{10}(3) = 3^3$.
\end{abstract}

\section{Introduction}

The $n$th $s$-gonal number, with $s \geq 3$, which we denote by $\mathcal{P}_s(n)$, is given by the formula \[ \mathcal{P}_s(n) = \frac{(s-2)n^2-(s-4)n}{2}. \] Polygonal numbers have been studied since antiquity \citep[pp.1--39]{history} and relations between different polygonal numbers and perfect powers have received much attention in the literature (see \citep{poly}, for example, and references therein). Kim, Park, and Pint{\'e}r \citep[Theorem 1.2]{poly} find all solutions to the equation $\mathcal{P}_s(n) = y^m$ when $m>2$ and  $s \in \{3,5,6,8,20 \}$ for integers $n$ and $y$. We extend this result (for $m>1$) to the case $s=10$, that of decagonal numbers. 

\begingroup
\renewcommand\thetheorem{1}
\begin{theorem}\label{mainthm}
All solutions to the equation \begin{equation}\label{maineq} \mathcal{P}_{10}(n) = y^m,\quad n,y,m \in \Z,~ m > 1, \end{equation} satisfy $n=y=0, n=\abs{y}=1,$ or $n=y=m=3$.

In particular, the only decagonal number $>1$ expressible as a perfect $m$th power with $m>1$ is $\mathcal{P}_{10}(3) = 3^3$.
\end{theorem}
\endgroup

We will prove Theorem \ref{mainthm} by carrying out a descent argument to obtain various ternary Diophantine equations, to which one may associate Frey elliptic curves. The difficulty in solving the equation $\mathcal{P}_{s}(n) = y^m$ for a fixed value of $s$ is due to the existence of the trivial solution $n=y=1$ (for any value of $m$). We note that adapting our method of proof also works for the cases $s \in \{3,5,6,8,20 \}$ mentioned above, but will not extend to any other values of $s$ (see Remark \ref{rem}).

\section{Descent and Small Values of $m$}

We note that it will be enough to prove Theorem \ref{mainthm} in the case $m=p$, prime. We write (\ref{maineq}) as \begin{equation}\label{eq2} n(4n-3) = y^p, \quad n,y \in \Z,~ p \text{ prime}, \end{equation} and suppose that $n,y \in \Z$ satisfy this equation with $n \ne 0$.

\subsubsection*{Case 1: $3 \nmid n$.}

If $3 \nmid n$, then $n$ and $4n-3$ are coprime, so there exist coprime integers $a$ and $b$ such that \[ n = a^p \quad \text{ and } \quad 4n-3 = b^p. \] It follows that \begin{equation} \label{case1} 4a^p-b^p = 3. \end{equation} If $p=2$ we see that $(2a-b)(2a+b)=3$,  so that $a = b = \pm 1$ and so $n=\abs{y}=1$. If $p=3$ or $p=5$, then using \texttt{Magma}'s \citep{magma} Thue equation solver, we find that $a=b=1$ also.

\subsubsection*{Case 2: $3 \parallel n$.}

Suppose $3 \parallel n$ (i.e. $\mathrm{ord}_3(n)=1$). Then after dividing (\ref{eq2}) by $3^{\mathrm{ord}_3(y)p}$, we see that there exist coprime integers $t$ and $u$, with $3 \nmid t$ such that \[n = 3t^p \quad \text{ and } \quad 4n-3 = 3^{p-1}u^p. \] Then \begin{equation} \label{case2} 4t^p-3^{p-2}u^p = 1. \end{equation} If $p=2$ we have $(2t-u)(2t+u)=1$ which has no solutions. If $p=3$, then $4t^3-3u^3=1$, and using \texttt{Magma}'s Thue equation solver, we verify that $u=t=1$ is the only solution to this equation, and this gives $n=y=3$. If $p=5$, we use \texttt{Magma}'s Thue equation solver to see that there are no solutions. 

\subsubsection*{Case 3: $3^2 \mid n$.}

If $3^2 \mid n$ then $3 \parallel 4n-3$, and arguing as in Case 2, there exist coprime integers $v$ and $w$, with $3 \nmid w$ such that \[n = 3^{p-1}v^p \quad \text{ and } \quad 4n-3 = 3w^p. \] So \begin{equation} \label{case3} 4 \cdot  3^{p-2}v^p- w^p = 1. \end{equation} If $p=2$, then as in Case 2 we obtain no solutions. If $p=3$ or $p=5$ then we use \texttt{Magma}'s Thue equation solver to verify that there are no solutions with $v \ne 0$.

\section{Frey Curves and the Modular Method}

To prove Theorem \ref{mainthm}, we will associate Frey curves to equations (\ref{case1}), (\ref{case2}), and (\ref{case3}), and level-lower to obtain a contradiction. We have considered the cases $p=2,3,$ and $5$ in Section 2, and so we will assume $m=p$ is prime, with $p \geq 7$.

We note that at this point we could directly apply \citep[Theorem 1.2]{pp2}  to conclude that the only solutions to (\ref{tern1}) are $a=b=1$, giving $n=1$, and apply \citep[Theorem 1.2]{consecutive} that (\ref{tern2}) and (\ref{tern3}) have no solutions. The computations for (\ref{tern1}) are not explicitly carried out in \citep{pp2}, so for the convenience of the reader and to highlight why the case $s=10$ is somewhat special, we provide some details of the arguments. 

\subsubsection*{Case 1: $3 \nmid n$.}

We write (\ref{case1}) as \begin{equation} \label{tern1} -b^p +4a^p = 3 \cdot 1^2, \end{equation} which we view as a generalised Fermat equation of signature $(p,p,2)$. We note that the three terms are integral and coprime. 

We suppose $ab \ne  \pm 1$. Following the recipes of \citep[pp.~26--31]{pp2}, we associate Frey curves to (\ref{tern1}). We first note that $b$ is odd, since $b^p = 4n-3$. If $a \equiv 1 \pmod{4}$, we set \[E_1: Y^2 = X^3 -3X^2+3a^pX. \] If $a \equiv 3 \pmod{4}$, we set  \[E_2: Y^2 = X^3 + 3X^2+3a^pX. \] If $a$ is even, we set \[E_3: Y^2+XY = X^3 - X^2 + \frac{3a^p}{16} X. \]

We level-lower each Frey curve and find that for $i=1,2,3,$ we have $E_i \sim_p f_i$, for $f_i$ a newform at level $N_{p_i}$, where $N_{p_1} = 36, N_{p_2} = 72,$ and $N_{p_3} = 18$. The notation $E \sim_p f$ means that the mod-$p$ Galois representation of $E$ arises from $f$. There are no newforms at level $18$ and so we focus on the curves $E_1$ and $E_2$. There is a unique newform, $f_1$, at level $36$, and a unique newform, $f_2$, at level $72$.

The newform $f_1$ has complex multiplication by the imaginary quadratic field $\Q(\sqrt{-3})$. This allows us to apply \citep[Proposition 4.6]{pp2}. Since $2 \nmid ab$ and $3 \nmid ab$, we conclude that $p=7$ or $13$, and that all elliptic curves of conductor $2p$ have positive rank over $\Q(\sqrt{-3})$. However, it is straightforward to check that this is not the case for $p=7$ and $13$. We conclude that $E_1 \not\sim_p f_1$.

Let $F_2$ denote the elliptic curve with Cremona label 72a2 whose isogeny class corresponds to $f_2$. This elliptic curve has full two-torsion over the rationals and has $j$-invariant $2^{4} \cdot 3^{-2} \cdot 13^{3}$. We apply \citep[Proposition 4.4]{pp2}, which uses an image of inertia argument, to obtain a contradiction in this case too.

\begin{remark}\label{rem} The trivial solution $a=b=1$ (or $n=y=1$) corresponds to the case $i=1$ above. The only reason we are able to discard the isomorphism $E_1 \sim_p f_1$ is because the newform $f_1$ has complex multiplication. The modular method would fail to eliminate the newform $f_1$ otherwise. For each value of $s$, we can associate to  (\ref{maineq}) generalised Fermat equations of signature $(p,p,2)$, $(p,p,3)$, and $(p,p,p)$. We found we could only obtain newforms with complex multiplication (when considering the case corresponding to the trivial solution) when $s = 3, 6, 8, 10,$ or $20$. A similar strategy of proof also works for $s=5$ using the work of Bennett \citep[p.~3]{approx} on equations of the form $(a+1)x^n-ay^n = 1$ to deal with the trivial solution.
\end{remark}

\subsubsection*{Case 2: $3 \parallel n$.}

We rewrite (\ref{case2}) as \begin{equation} \label{tern2} 4t^p-3^{p-2}u^p = 1 \cdot 1^3 , \end{equation} which we view as a generalised Fermat equation of signature $(p,p,3)$. The three terms are integral and coprime. We suppose $tu \ne \pm 1$. Using the recipes of \citep[pp.~1401--1406]{pp3}, we associate to (\ref{tern2}) the Frey curve \[E_4: Y^2 + 3XY - 3^{p-2}u^p \, Y    = X^3.   \] We level-lower $E_4$ and find that $E_4 \sim_p f$, where $f$ is a newform at level $6$, an immediate contradiction, as there are no newforms at level $6$.

\subsubsection*{Case 3: $3^2 \mid n$.}

We rewrite (\ref{case3}) as \begin{equation} \label{tern3} -w^p + 4 \cdot  3^{p-2}v^p = 1 \cdot 1^3, \end{equation} which we view as a generalised Fermat equation of signature $(p,p,3)$. The three terms are integral and coprime. We suppose $vw \ne \pm 1$. The Frey curve we attach to (\ref{tern3}) is 
\[ E_5: Y^2 + 3XY + 4 \cdot 3^{p-2} v^p \, Y = X^3. \] We level-lower, and find that $E_5 \sim_p f$, where $f$ is a newform at level $6$, a contradiction as in Case 2.

This completes the proof of Theorem \ref{mainthm}.

\bibliographystyle{plainnat}

\begin{thebibliography}{26}
\providecommand{\natexlab}[1]{#1}
\providecommand{\url}[1]{\texttt{#1}}
\expandafter\ifx\csname urlstyle\endcsname\relax
  \providecommand{\doi}[1]{doi: #1}\else
  \providecommand{\doi}{doi: \begingroup \urlstyle{rm}\Url}\fi

\bibitem{approx}
M.~Bennett. 
\newblock Rational approximation to algebraic numbers of small height: the {D}iophantine equation {$ \, \abs{ax^n+by^n} = 1$}.
\newblock \emph{J. Reine Angew. Math.}, 535:\penalty0 1--49, 2001.

\bibitem{consecutive}
M.~Bennett.
\newblock Products of consecutive integers.
\newblock \emph{Bull. Lond. Math. Soc.}, 36\penalty0(5):\penalty0  683--694, 2004.

\bibitem{pp2}
M.~Bennett and C.~Skinner.
\newblock Ternary {D}iophantine equations via {G}alois representations and modular forms.
\newblock \emph{Canad. J. Math.}, 56\penalty0 (1):\penalty0 23--54, 2004.

\bibitem{pp3}
M.~Bennett, V.~Vatsal, and S.~Yazdani. 
\newblock Ternary {D}iophantine equations of signature {$(p, p, 3)$}.
\newblock \emph{Compos. Math.}, 140\penalty0(6):\penalty0 1399--1416, 2004.


\bibitem[Bosma et~al.(1997)Bosma, Cannon, and Playoust]{magma}
W.~Bosma, J.~Cannon, and C.~Playoust.
\newblock The {M}agma algebra system. {I}. {T}he user language.
\newblock \emph{J. Symbolic Comput.}, 24\penalty0 (3-4):\penalty0 235--265, 1997.

\bibitem{history}
L.~Dickson
\newblock \emph{History of the Theory of Numbers. Volume II: Diophantine Analysis}.
\newblock Dover Publications, New York, 2005.

\bibitem{poly}
D.~Kim, K.~Park, and A.~Pint{\'e}r.
\newblock A {D}iophantine problem concerning polygonal numbers.
\newblock \emph{Bull. Aust. Math. Soc.}, 88\penalty0 (2):\penalty0 345--350, 2013.
\end{thebibliography}

\Addresses

\end{document}